\documentclass[letterpaper, 10 pt, conference]{ieeeconf}  

\usepackage{epsfig}
\usepackage{subcaption}
\usepackage{calc}
\usepackage{amssymb}
\usepackage{amstext}
\usepackage{amsmath}

\usepackage{amsthm}
\usepackage{multicol}
\usepackage{pslatex}
\usepackage{float}
\usepackage{url}
\usepackage{hyperref}

\usepackage{algorithm}
\usepackage{algpseudocode}

\theoremstyle{definition}
\newtheorem{definition}{Definition}

\IEEEoverridecommandlockouts                              
\title{\LARGE \bf
	Clustering as a means of leader selection in consensus networks
}

\author{Natalia Basimova$^{1}$ and  Pavel Chebotarev$^{2}$
	\thanks{$^{1}$Author is with Skolkovo Institute of Science and Technology,
        30 Bolshoi Boulevard, Bld.1, Moscow 121205, Russia
		{\tt\small basimova.nf@phystech.edu}}%
	\thanks{$^{2}$Author with Institute of Control Sciences of the Russian Academy of Sciences,
		65 Profsoyuznaya str., Moscow 117997, Russia
		{\tt\small upi@ipu.ru}}%
}

\begin{document}

	
	 \maketitle
      \pagestyle{plain}
      \thispagestyle{plain}

	
	 \begin{abstract}
	 	In the leader-follower approach, one or more agents are selected as leaders who do not change their states or have autonomous dynamics and can influence  other agents, while the other agents, called followers, perform a simple protocol based on the states of their neighbors. This approach provides a natural link between control theory and networked agents with their input data. Despite the fact that the leader-follower approach is widely used, the fundamental question still remains: how to choose leaders from a set of agent. This question is called the problem of choosing leaders. There is still no selection algorithm that is both optimal under a natural criterion and fast. In this paper, for agents that obey a linear consensus protocol, we propose to choose leaders using graph nodes'  clustering algorithms and show that this method is the most accurate among the fast existing algorithms of choosing leaders.
	\end{abstract}	

	
	
\section{\uppercase{Introduction}}
\label{sec:introduction}
	
	A multi-agent system is a system formed by a number of interacting agents. Multi-agent systems can be used to solve problems that are difficult or impossible to solve with a single agent or a monolithic system. Examples of such tasks are emergency response \cite{1}, modeling of social structures \cite{2}, and many others.
	
	Movement of several agents to a given control point while maintaining a relationship between their states is an important problem in robotics. For example, in applications, including spaceships, unmanned aerial vehicles (UAVs), and mobile robots, agents often need to move to a landmark or target point while maintaining the shape of their formation.
	
	The leader-follower approach is one of the ways to solve this type of formation control problems \cite{3}.
	
	Recently, several methods of solving the problems of choosing leaders in networked multi-agent systems have been studied, in particular, within the framework of linear consensus protocols.
	Some graph indices, such as node degree, have been tested in terms of leader performance \cite{4}. The problem of choosing leaders can be solved either ascendingly (choosing one leader and adding other leaders until the set of leaders becomes optimal) \cite{main-paper}, or descendingly: the optimal set of leaders is sought by excluding agents from the entire set of agents \cite{6}. To control UAVs, an ``online'' leader change was applied in order to improve the rate of convergence of agents' positions \cite{5}.
	
	The task of selecting leaders from the set of agents is similar to a node clustering task in a graph if each leader  is treated as an agent responsible for its own cluster \cite{7}. Therefore, it is reasonable to compare existing leader selection algorithms with clustering methods combined with choosing the ``main'' node in each cluster. The most popular clustering methods are spectral clustering based on optimization of graph slices \cite{8} and the $k$-means method, which seeks to minimize the total quadratic deviation of cluster points from the centers of these clusters \cite{9}.
	
	In this paper, we consider the problem of choosing leaders in the case where they preserve their initial states, i.e., they are static. The problem of  selecting mobile leaders was considered, for example, in \cite{10}.
	
	The work consists of two main parts. The first one is a simulation in ideal conditions and the second one is a simulation in real conditions using a model of e-puck2 mobile robots \cite{e-puck2}. ``Real conditions'' means a system with some restrictions, for example, with a limited possible speed of the robots.
	
	The aim of the paper is to compare existing algorithms for choosing leaders with the method based on graph nodes' clustering with respect to accuracy and convergence time in ideal and real conditions.
	
	\section{\uppercase{Model}}
	
	We will consider a multi-agent system with a corresponding connection graph $\mathcal{G} = (\mathcal{V}, \mathcal{E})$. The graph $\mathcal{G}$ is directed and represents the system structure. We will identify the node $ v_i $, $ i = 1, ..., V $, with the $i$-th agent. Agent $i$ depends on agent $j$ if there is a directed edge $e = (v_i, v_j) \in \mathcal{E}$. Each edge $e \in \mathcal{E}$ has its  weight $w$, thereby we have a weighted adjacency matrix $A$ of the system \cite{big-book}.
	
	Let matrix $L$ be the Laplacian matrix \cite{lapl} of the graph $\mathcal{G}$ with adjacency matrix $A$. Each agent $i$ has its state $x_i$. In this paper, $x_i$ denotes agent's coordinate. The linear consensus protocol has the form
	\begin{equation}\label{ed:linear-consesus-equation}
	\bold{\dot{x}} = - L \bold{x},
	\end{equation}
	where $\bold{x}=(x_1,\ldots, x_V)^T$ is the state vector of the system.
	
	Leaders represent stubborn agents: they get their states initially, and these states do not change over time. Respectively, we will divide the general state of the system into the state of leaders and the state of followers: $\bold{x}=  \Big[
	\begin{array}{cc}
	\bold{x_F} \\
	\bold{x_L} \\
	\end{array}
	\Big],$ where $\bold{x_F}$ is the state of followers and $\bold{x_L}$ is the state of leaders. Thus, the linear consensus protocol (\ref{ed:linear-consesus-equation}) can be rewritten in the following way:
	
	\begin{equation}
	\Big[
	\begin{array}{cc}
	\bold{\dot{x}_F} \\
	\bold{\dot{x}_L} \\
	\end{array}
	\Big]
	=
	-\left( \begin{array}{cc}
	L_{\text{FF}} & L_{\text{FL}}  \\
	L_{\text{LF}} & L_{\text{LL}}  \\
	\end{array} \right)
	\Big[
	\begin{array}{cc}
	\bold{x_F} \\
	\bold{x_L} \\
	\end{array}
	\Big].
	\end{equation}
	
	Leaders preserve their states, hence, $L_{\text{LL}}$ and $L_{\text{LF}}$ are zero, so that $ \bold{\dot{x}_L} $ is zero and $ \bold{x_L} $ is constant. Then, $L_{\text{FF}}$ is the grounded Laplacian of the system \cite{grounded}. In this paper, $L_{\text{FF}}$ is assumed to be symmetric. Thus, the followers influence each other equally.
	The convergence rate of the system can be asymptotically measured by the smallest eigenvalue of the grounded Laplacian $L_{\text{FF}}$ \cite{main-paper}, therefore, further in the paper we will identify these quantities.
	
	
	
	Next, for simulation purposes, we define convergence time as follows:
	\begin{definition}\label{def:time}
		Let $\bold{x^*}$ be a limit state of the system. The system reaches $\bold{x^*}$  with fixed error $e$ in convergence time $t_e$, if
		\begin{equation}
		||\bold{x}(t_+) - \bold{x^*}||_2 \leq e \quad \forall t_+ \geq t_e,
		\end{equation}
		where $||\bold{y}||_2$ is the Euclidean norm of vector $ \bold{y}$.
	\end{definition}

	\section{\uppercase{Algorithms}}\label{sec-algorithms}
	We will compare six algorithms, also we will study the algorithmic complexity of them, after that we will consider two additional algorithms presented in \cite{algorithmkn3, algorithmn4}. The algorithm in \cite{algorithmkn3} is based on theoretical calculations, so it is optimal, however its algorithmic complexity is higher than that of the best realization of the greedy algorithm in \cite{algorithmn4}, whose complexity in turn is higher than that of the \textit{k-means} algorithm.
	
	These six algorithms are:
	\begin{enumerate}
		\item \textit{Choose k-leader algorithm}. A pseudocode description of \textit{Choose k-leader algorithm} is given in \cite{main-paper}. Its algorithmic complexity is between $ O(kn^3) $ and $ O(kn^4) $ due to the need to calculate the eigenvalues of the Laplacian matrix and then choose the smallest eigenvalue $k$ times.
		\item \textit{Random leader selection}. \textit{k} numbers are randomly selected  without repetitions using uniform distribution from $ 1 $ to $ V $; the corresponding agents are chosen as leaders. The algorithmic complexity is $ O(n) $.
		\item \textit{Leaders with max degree}. \textit{k} nodes with maximum indegrees are chosen from the connection graph. The algorithmic complexity is $ O(n^2) $, because we need to calculate the degree of each node.
		\item \textit{Leaders with average degree}. \textit{k} nodes
		 with average indegrees are chosen from the connection graph. The algorithmic complexity is $ O(n^2) $ as well.
		\item \textit{k-means algorithm}. The connection graph is split into \textit{k} clusters with $k$-means algorithm. Then, the closest node to the center of the cluster is chosen as the leader of this cluster. If this node is already chosen as leader for another cluster, then the second closest node is selected and so on. The $k$-means clustering is taken from \texttt{Python} library \texttt{sklearn} with default parameters; its average algorithmic complexity is $ O(kn)$, the worst complexity is $ O(n^{k+1})$ \cite{sklearn}. However, the worst complexity is not reached because of the fixed maximum number of iterations.
		\item \textit{Huge random selection}. Firstly, 10000 numbers from 1 to $ {n\choose k} $ are uniformly chosen and each number is an index in an imaginary list in which all possible selections of $k$ leaders are arranged in lexicographic order; knowing the index, a selection can be made. Thus, 10000 selections are chosen. The algorithm iterates over all selected sets of leaders and finds the set with the maximum smallest eigenvalues of the corresponding grounded Laplacian $L_{\text{FF}}$. If the number of different selections is not uniformly limited, then the algorithm has combinatorial complexity. Therefore, the upper bound of  10000 was taken in order to reduce the execution time to a reasonable one. Obviously, the algorithm does not guarantee the best result, but as we will see in the experiments, its results are better than those of the previous algorithms.
	\end{enumerate}
	
	As for the remaining two algorithms, greedy algorithm in \cite{algorithmn4} with the best realization has an algorithmic complexity of $O(n^3)$ and the complexity without these improvements in realization is $O(n^4)$.
	The exact algorithm in \cite{algorithmkn3} for the $k$-leader selection problem has a polynomial algorithmic complexity of $O(n^3k)$ in the case of path or ring graphs taken.
	
	Thus, it can be seen that the \textit{k-means algorithm} has the smallest average asymptotic complexity among non-random algorithms.
	
	
	\section{\uppercase{Simulation}}
	
	Consider a multi-agent system containing 100 agents uniformly distributed in a $ 10 \times 10 $ m$^2$ square. Let two agents-followers be symmetrically connected if and only if the initial distance between them is no more than $ 3 $ m. The weights of these connections are chosen uniformly from the interval (0,~50]. The weighted adjacency matrix generated in this way does not change later.
	
	In all experiments, we measured the convergence time (see Definition \ref{def:time}). Since agents had $x$ and $y$ coordinates, the total deviations along both axes from the limiting state were calculated, and when both of these deviations do not exceed the fixed error, we admit that the system reached the limit. Based on the obtained convergence times, we conclude that one algorithm is better than another.
	Additionly, we measure the convergence rates for all algorithms, which are identified with the smallest eigenvalues of the corresponding grounded Laplacian matrices $L_{\text{FF}}$; the results are provided below.
	
	Each experiment was repeated 30 to 100 times depending on the number of leaders $k$ (from $1$ to $90$); the obtained convergence time was averaged for each $k$.
	
	To bring the system closer to reality, physical limitations were imposed: these were the limits on the maximum speed. Simulation was performed in Webots, an open source robot simulator, which contains various robot models, including a model of e-puck2 robots. The e-puck2 model was used for simulation with physical constraints that the maximum speed is limited to $15.4$ cm/s.
	
	Experiments with physical limitations were carried out in Webots on \texttt{C} language, without these restrictions on \texttt{Python}. More details about the experiments carried out can be found in the Appendix.
	
	\subsection{Experiments' settings}
	
	For all experiments, the convergence error (see Definition~\ref{def:time}) was set to be $ 5 \times 10 ^ {- 8} $ cm. This error was chosen experimentally, as it allows to achieve fairly accurate results within a reasonable time.
	
	The number of leaders ranged from $ 1 $ to $ 90 $, but the charts show values only from $ 1 $ to $ 9 $, since the results for larger numbers of leaders are similar. Convergence time slowly discreases and convergence rate slowly increases in the number of leaders for all dependencies.
	
	The simulation time step for both experiments with physical limitations and without them was equal to 1 ms, which is the minimum available time step in Webots.

	\subsection{Convergence time of the algorithms}\label{main-exp-sec}
	
	In this section, we present the results on the convergence time depending on the number of leaders either with constraints on speed or without them.

\begin{figure}[thpb]
	\centering
	\begin{subfigure}{.5\textwidth}
		\centering
		\includegraphics[width=\linewidth]{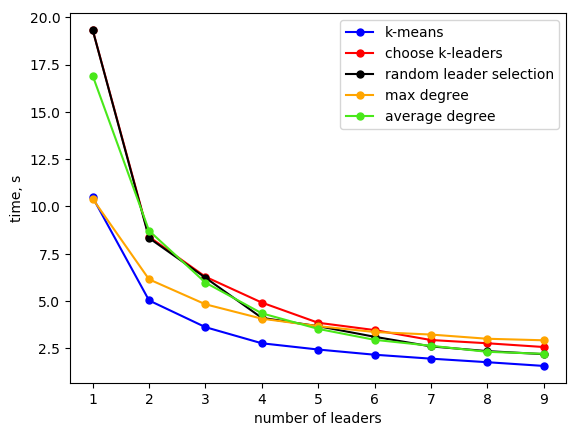}
		\caption{}
		\label{fig:sub1}
	\end{subfigure}%
	\hfill
	\begin{subfigure}{.5\textwidth}
		\centering
		\includegraphics[width=\linewidth]{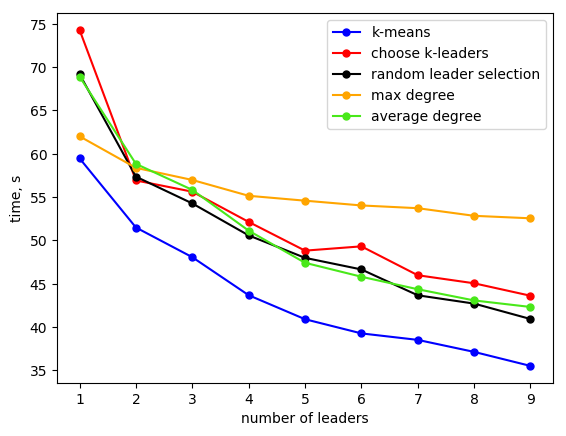}
		\caption{}
		\label{fig:sub2}
	\end{subfigure}
	\caption{Comparison of \textit{k-means} with four other algorithms. (a)  Experiments without physical limitations. (b) Experiments with physical limitations on speed.}
	\label{fig:fisrt-exp}
\end{figure}

First, we compare algorithms 1 to 5 because their complexities are much smaller than the complexity of \textit{huge random selection} algorithm that is defined by constant number 10000. As Fig. \ref{fig:fisrt-exp} shows, \textit{k-means algorithm} gives the best result and \textit{random leader selection}, \textit{choose k-leader algorithm}, and \textit{leaders with average degree} give approximately similar results. Finally, \textit{leaders with max degree} gives the worst result in the experiments with physical restrictions. However, the algorithm shows good results with a small number of leaders in the experiment without physical restrictions. These results should be investigated separately. Note that the difference between the algorithms is clearly seen in the experiments with physical limitations, however the convergence rate, which in this study is identified with the smallest eigenvalue of the grounded Laplacian, is the same in both experiments, therefore the convergence rate should not be the only decisive factor when choosing an algorithm.

Also the convergence time of experiments with physical limitations (see Fig. \ref{fig:sub2}) does not smoothly decrease with increasing number of leaders; this may indicate that more averaging is needed.

\begin{figure}[thpb]
	\centering
	\begin{subfigure}{.5\textwidth}
		\centering
		\includegraphics[width=1\linewidth]{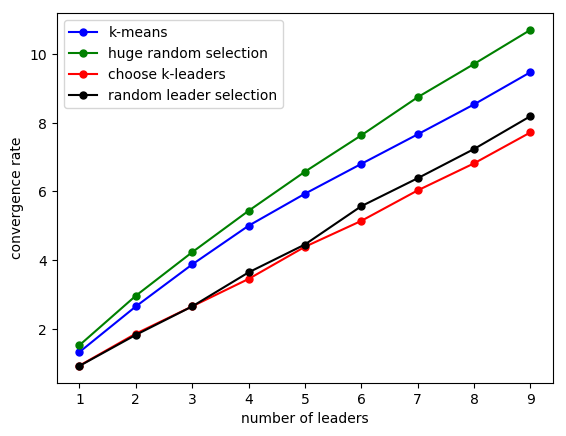}
		\caption{}
		\label{fig:sub3}
	\end{subfigure}%
\hfill
	\begin{subfigure}{.5\textwidth}
		\centering
		\includegraphics[width=1\linewidth]{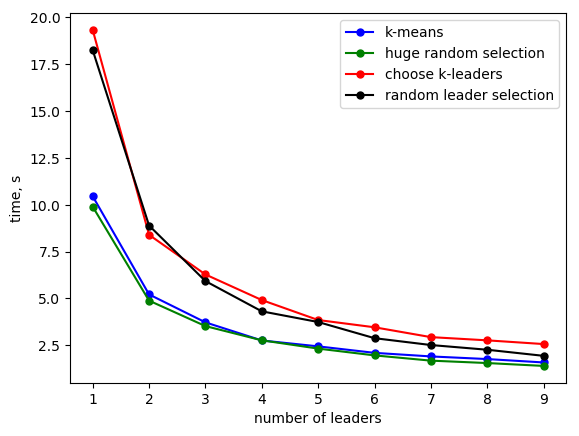}
		\caption{}
		\label{fig:sub4}
	\end{subfigure}
\hfill
    \begin{subfigure}{.5\textwidth}
		\centering
		\includegraphics[width=1\linewidth]{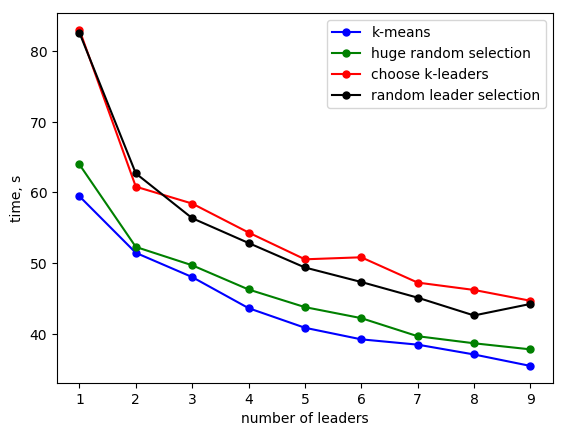}
		\caption{}
		\label{fig:sub5}
    \end{subfigure}
	\caption{Comparison of \textit{huge random selection}  with three previous algorithms. (a) Convergence rates in the experiments. (b)  Time in the experiment without physical limitations. (c) Time in the experiment with physical limitations on speed.}
	\label{fig:second-exp}
\end{figure}

Second, we compare \textit{k-means}, \textit{choose k-leader} and \textit{random leader selection} algorithms with \textit{huge random selection}. We selected these algorithms because they gave the best or approximately similar results in the previous experiments and we do not compare all the algorithms so as not to clutter up the charts. As we can see in Fig.  \ref{fig:sub3}, \textit{huge random selection} gives the best convergence rate and \textit{k-means}' result is quite close to it. Also we can see from the simulation that these two algorithms give the best results concerning the convergence time, though \textit{k-means algorithm} is better in the experiments with physical restrictions. This is apparently related to the stop simulation criterion which depends on the convergence error (see Definition \ref{def:time}).

\begin{figure}[thpb]
	\centering
	\begin{subfigure}{.5\textwidth}
		\centering
		\includegraphics[width=\linewidth]{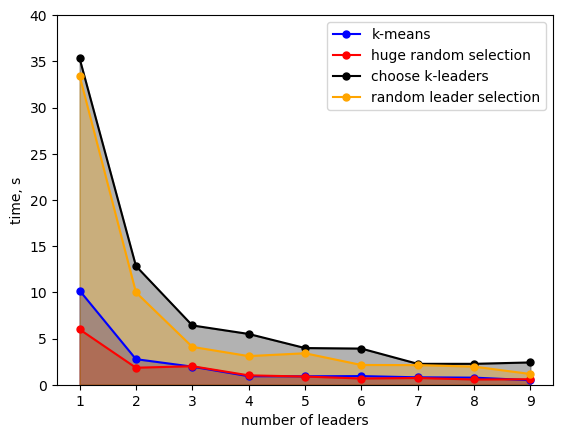}
		\caption{}
		\label{fig:sub6}
	\end{subfigure}%
	\hfill
	\begin{subfigure}{.5\textwidth}
		\centering
		\includegraphics[width=\linewidth]{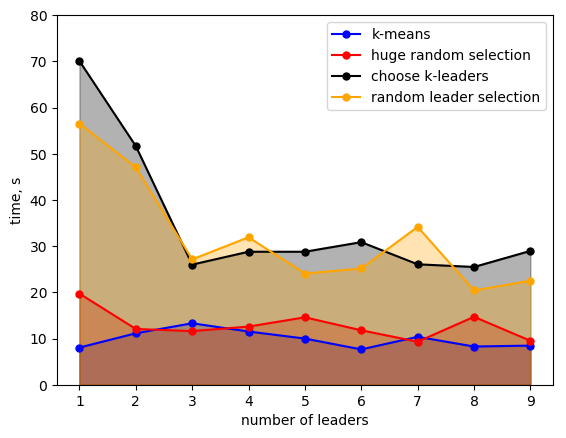}
		\caption{}
		\label{fig:sub7}
	\end{subfigure}
	\caption{Scatter of convergence time. (a) The experiment without physical limitations. (b) The experiment with physical limitations}
	\label{fig:last-exp}
\end{figure}

Finally, let us look at the scatter of the received data, which will give an understanding of the stability of the algorithm for real use. On Fig. \ref{fig:last-exp}, the difference between the maximum and the minimum convergence time is plotted along the time axis. It can be seen that the ranges of time for \textit{k-means} and \textit{huge random selection} algorithms are much smaller than those of  \textit{choose k-leader} and \textit{random leader selection} algorithms; for the remaining algorithms the range is quite similar to the last two.

\section{\uppercase{Conclusion}}
In this paper, we compared several algorithms for solving the leader selection problem with regard to the leader-follower consensus protocol for both systems with physical constraints and idealized ones. The following algorithms were compared:
\begin{itemize}
	\item The \textit{$k$-leader selection} algorithm  \cite{main-paper},
	\item the \textit{random choice},
	\item the \textit{choice of agents with the maximum indegrees} in the agents' dependency graph,
	\item the \textit{choice of the agents with average indegrees} in the agents' dependency graph,
	\item the \textit{$k$-means algorithm}, with the subsequent selection of nodes that are closest to the center of each cluster.
\end{itemize}

Additionally, all algorithms were compared with \textit{huge random selection} algorithm. It turned out that the clustering algorithm showed the best results in all experiments, except for \textit{huge random selection} algorithm whose algorithmic complexity is determined by the upper bound on the number of different selections, without this bound, the complexity is combinatorial.

In future work, it would be interesting to apply other clustering algorithms, as well as an alternative version of $k$-means algorithm, where the leader in each cluster is chosen as the node with the maximum degree in the cluster.

In preliminary experiments, we observed that the agglomerative clustering algorithm gave similar results to those of the $k$-means algorithm, though the algorithm complexity of the agglomerative clustering is higher than the complexity of the $k$-means algorithm \cite{clusters}. Thus, the main conclusion is that clustering algorithms are useful in solving the leader selection problem. However, we also tested the spectral clustering algorithm \cite{spectral} and it gave results much worse than those of \textit{leaders with max degree} algorithm.


In this work, we did not consider systems with mobile leaders or dynamic connections, systems with acceleration, or systems with shape preserving. The problems of leader selection for them are in many respects similar to the simplest problems studied in this paper, however, of course, they are more complicated both theoretically and experimentally. In a further development of this work, such problems  can be analyzed by applying clustering algorithms.

\section*{APPENDIX}
Simulation on \texttt{Python} and Webots consisted of 5 main steps:
\begin{enumerate}
	\item Generation of coordinates;
	\item Generation of the adjacency matrix based on these coordinates and the connection condition;
	\item Generation of the Laplacian matrix;
	\item A transformation of this Laplacian matrix according to the type of an algorithm (see Section \ref{sec-algorithms});
	\item Run of the simulation.
\end{enumerate}

The first four steps were the same for both types of the experiments and were implemented on \texttt{Python}.

After the generation of the coordinates and the Laplacian matrix the simulation started in a huge cycle, which can contain a maximum of $N$ iterations; in our experiments, $N = 2 \times 10^6$.

Firstly the final limit state $\bold{x^*}$ was calculated simply as $e^{-LNt_s}\bold{x}(0)$, so if the system did not reach this limit state in $N$ steps or less, new coordinates and the corresponding Laplacian matrix were generated and simulation started again.

\subsection*{Experiments  without physical limitations}
On each step of the simulation on \texttt{Python}, agents' coordinates were updated as:
$$
\bold{x} = \bold{x}  - L \bold{x} t_s,
$$
where $t_s$ is a simulation time step and equals to 1 ms in our experiments. After this update, the distance between the current state and the limit state was calculated according to Definition~\ref{def:time} and compared with the fixed error. If this distance was less than the fixed error, the simulation stoped, otherwise the new iteration began.

\subsection*{Experiments  with physical limitations}
Simulation on Webots was similar to the one on \texttt{Python} with the difference that before the each simulation step $t_s$, agents updated their new speeds as $\bold{\dot{x}} = \text{min}(- L \bold{x}, \text{15.4 cm/s}), $ with 15.4 cm/s being the maximum speed. Webots updates the speed with the last speed, acceleration and motor torque taken into account. After a simulation time step, a predeterminated agent (``inspector'') obtained the other agents' coordinates, calculated the distance  from the limit state, compared it with the fixed error and decided if simulation required continuation. If so, the inspector recalculated the speeds of all agents, sent them back, and simulation step repeated again.

Due to the way of the simulation and the implementation, no time was wasted on exchanging information, this time must be taken into account when implementing the algorithm on real robots.




\end{document}